\documentclass[12pt]{amsart}
\usepackage{amssymb}
\usepackage[T1]{fontenc}
\usepackage[all]{xy}
\UseRawInputEncoding
\numberwithin{equation}{section}

\newcommand{\Q}{{\mathbb {Q}}}

\newcommand{\G}{{\bf{G}}}
\newcommand{\T}{{\bf{T}}}

\newcommand{\R}{{\mathbb{R}}}
\newcommand{\Z}{{\mathbb{Z}}}
\newcommand{\C}{{\mathbb{C}}}

\newcommand{\N}{{\mathbb{N}}}

\newcommand{\OO}{{\mathcal O}}

\theoremstyle{plain}
\newtheorem{thm}{Theorem}[section]
\newtheorem{lem}[thm]{Lemma}
\newtheorem{prop}[thm]{Proposition}
\newtheorem{cor}[thm]{Corollary}

\theoremstyle{definition}

\title[Norm forms]{Action of maximal tori on homogeneous spaces with a view to number theory}

\title[Norm forms]{Characterization of norm and quasi-norm forms in $S$-adic setting}

\author{George Tomanov}
\address{Institut Camille Jordan, Universit\'e Claude Bernard - Lyon
I, B\^atiment de Math\'ematiques, 43, Bld.
 du 11 Novembre 1918,
69622 Villeurbanne Cedex, France {\tt tomanov@math.univ-lyon1.fr}}

\begin{document}

\maketitle

\begin{abstract} The goal of the present paper is to characterize the norm and quasi-norm forms defined over an arbitrary number field $F$ in terms of their values at the $S$-integer points, where $S$ is a finite set of valuations of $F$ containing the archimedean ones. In this way we generalize the main result of the recent paper \cite{Toma4}, where the notion of a quasi-norm form is introduced when $F = \Q$ and $S$ is a singleton.
In complement, we exhibit some relations with problems and results in this area of research.
\end{abstract}

\section{Introduction} \label{Introduction}

The investigation of different types of (algebraic) forms via their
values at the integer points has a long history. Impressive results in this direction were obtained during the last decades thanks to reformulation of longstanding problems and conjectures in homogeneous dynamical terms. This approach allows to apply, in parallel with deep methods from the algebra and the algebraic geometry, powerful tools from the analysis and the ergodic theory. As an example, recall the pathbreaking G.Margulis' proof in 1986 (\cite{Margulis+Oppenheim} and \cite{Oppenheim survey}) of the A.Oppenheim conjecture formulated in 1929.
The result asserts that if $q$ is a quadratic form with real coefficients in $n \geq 3$ variables which is  nondegenerate, indefinite, and not multiple to a form with rational coefficients then for any $\varepsilon > 0$ there exists $\vec{z} \in \Z^n \setminus \{\vec{0}\}$ such that $|q(\vec{z})| < \varepsilon$. Its homogeneous dynamics reformulation (proved by Margulis) states that the closure of every orbit for the action by left translation of the orthogonal group $\mathrm{SO}(q)$ on the space of unimodular lattices $\mathrm{SL}_{n}(\R)/\mathrm{SL}_{n}(\Z)$ is an orbit of a subgroup of $\mathrm{SL}_{n}(\R)$ containing $\mathrm{SO}(q)$. It is known that the Oppenheim conjecture was motivated by Meyer's theorem that if $q$ is a quadratic form with coefficients from $\Q$ in $n \geq 5$ variables which is nondegenerate and indefinite
then $q$ represents zero over $\Q$ non-trivially, that is, $q(\vec{z}) = 0$ for some $\vec{z} \in \Q^n \setminus \{\vec{0}\}$. The $S$-adic version of Oppenheim's conjecture is proved by A.Borel and Gopal Prasad in \cite{BP}. The fact that when $q$ is indefinite and $n \geq 3$ the group $\mathrm{SO}(q)$ is generated by its $1$-parameter unipotent subgroups is crucial for the proof. In general, if $G$ is a connected Lie group, $\Gamma$ a lattice in $G$, and $H$ a connected subgroup of
$G$ generated by its $\mathrm{Ad}$-unipotent $1$-parameter subgroups, it is proved by M.Ratner \cite{Ratner} that the closure of every $H$-orbit on $G/\Gamma$ is \textit{homogeneous}, that is, it is an orbit itself of a closed connected subgroup of $G$ containing $H$. The result confirms a conjecture of M.S.Raghunathan. Its proof is based on
the so-called measure rigidity of unipotent flows on homogeneous spaces conjectured by S.G.Dani around 1980 and proved
 by M.Ratner \cite{Ratner2}.
Inspired by \cite{BP}, $S$-adic versions of \cite{Ratner} and \cite{Ratner2} are proved in \cite{MaT1}, \cite{MaT2},  \cite{Ratner3}, and \cite{Toma1}.

One of the examples which shows that the assumption $n \geq 3$ in the formulation of the Oppenheim conjecture is essential and can not be omitted is provided by the form $q(x_1,x_2) = x_1^2 - \alpha^2x_2^2$, where $\alpha$ is a badly approximable real number such that $\alpha^2 \notin \Q$.
(See \cite[1.2]{Margulis problems survey} and \cite[ch.1, \S5]{Schmidt}.) In this case $\mathrm{SO}(q)$ is $\R$-\textit{split} (that is, diagonalizable over $\R$) and the closure of the $\mathrm{SO}(q)$-orbit of the lattice $\Z^2$ in the space
$\mathrm{SL}_{2}(\R)/\mathrm{SL}_{2}(\Z)$ of unimodular lattices is compact and not homogeneous.
In general, if $G = \G(\R)$, where $\G$ is a linear algebraic group defined over $\R$, $\Gamma$ an irreducible lattice in $G$, and $T$ an $\R$-split algebraic sub-torus of $G$ such that $G/\Gamma$ does not admit, in a natural sense, rank $1$ $T$-invariant factors, it was believed for a certain time that, as in the case of unipotent flows treated in \cite{Ratner} and \cite{Ratner2}, the closure of every $T$-orbit on $G/\Gamma$ should be homogeneous \cite[Conjecture 1]{Margulis problems survey}. It turned out that this is not always the case. There are counter-examples in the following cases: for the action of $n-2$-dimensional tori on $\mathrm{SL}_{n}(\R)/\mathrm{SL}_{n}(\Z)$ with $n \geq 6$ in \cite{Mau}, for the action of maximal tori on $\mathrm{SL}_{3}(\R)/\mathrm{SL}_{3}(\Z)$ in \cite{Sca},  for the action of maximal tori $T = T_1 \times T_2$ on $\mathrm{SL}_{2}(\R) \times \mathrm{SL}_{2}(\R)/\Gamma$ in \cite{Toma1+}, and, finally, for arbitrary semisimple groups in general $S$-adic setting in \cite{Toma3}. \textit{The common point in all these results is that the $T$-orbits under consideration are not relatively compact.} For relatively compact orbits, \cite[Conjecture 1]{Margulis problems survey} is open and  represents a generalization of the 1955 conjecture of Cassels and Swinnerton-Dyer  \cite{CS} that the closure of every relatively compact orbit for the action of the diagonal subgroup of $\mathrm{SL}_{n}(\R)$, $n \geq 3$, on $\mathrm{SL}_{n}(\R)/\mathrm{SL}_{n}(\Z)$ is homogeneous. It is worth mentioning that the conjecture of Cassels and Swinnerton-Dyer implies the famous conjecture of Littlewood stated around 1930. (See \cite{CS}, \cite[Conjectures 7 and 8]{Margulis problems survey} and \cite[\S2]{Oppenheim survey}.) An important progress on Littlewood's conjecture is due to
M.Einsiedler, A.Katok and E.Lindenstrauss who proved  that the set of exceptions to Littlewood's conjecture has Hausdorff dimension zero \cite{EKL}.

The results of this paper concern the relatively compact orbits for the action of arbitrary (not always spit) maximal algebraic tori on $S$-adic homogeneous spaces.
We get a description of the norm and quasi-norm forms over \textit{arbitrary} number fields in terms of their values at the $S$-integer points.
 In this way, we generalize the results from \cite{Toma4} about norm and quasi-norm forms over the field of rational numbers.
 In comparison with \cite{Toma4} our arguments are more general and in certain circumstances simplified.

Further on, $F$ is a number field, $S$ is a finite set of valuations of $F$ containing the archimedean ones and $\mathcal{O}_S$ is the ring of $S$-integers of $F$. For every $v \in S$, let $F_v$ be the completion of $F$ with respect to $v$ and $\overline{F}_v$ be the algebraic closure of $F_v$. Denote $F_S = \underset{v \in S}{\prod}{F}_v$ and $\overline{F}_S = \underset{v \in S}{\prod}\overline{F}_v$. The field $F$ is diagonally embedded in $F_S$. Let $F_S[\vec{x}]$ (respectively, $\overline{F}_S[\vec{x}]$) be the ring of polynomials in $m \geq 2$ variables $\vec{x} = (x_1, \dots, x_m)$. Then $F_S[\vec{x}] = \underset{v \in S}{\prod} F_v[\vec{x}]$ (respectively, $\overline{F}_S[\vec{x}] = \underset{v \in S}{\prod} \overline{F}_v[\vec{x}]$). The ring $F[\vec{x}]$ is diagonally embedded in $F_S[\vec{x}]$.
If $p(\vec{x}) \in F_S[\vec{x}]$ then $p(\vec{x}) = (p_v(\vec{x}))_{v \in S}$, where $p_v(\vec{x}) \in F_v[\vec{x}]$.
  By an \textit{$S$-form of degree $n$} (or, simply, a form of degree $n$ when $"S"$ is implicit)
 we mean a homogeneous polynomial $f(\vec{x}) \in F_S[\vec{x}]$ which factors as $f(\vec{x})= \overset{n}{\underset{i = 1}{\prod}}l_i(\vec{x})$, where $l_1(\vec{x}), \dots, l_n(\vec{x})$ are linearly independent over $\overline{F}_S$ linear forms in $\overline{F}_S[\vec{x}]$. Hence, if $f(\vec{x}) = (f_v(\vec{x}))_{v \in S}$ then every  $f_v(\vec{x})$ is an $F_v$-\textit{form} decomposable over $\overline{F}_v$, that is, $f_v(\vec{x}) = \overset{n}{\underset{i = 1}{\prod}}l_{v,i}(\vec{x})$, where $l_{v,1}(\vec{x}), \dots, l_{v,n}(\vec{x})$ are linearly independent linear forms in $\overline{F}_v[\vec{x}]$.

An $S$-form $f(\vec{x}) = (f_v(\vec{x}))_{v \in S}$ is \textit{decomposable over} $F_S$ if $l_{v,i}(\vec{x}) \in F_v[\vec{x}]$ for all $v$ and $i$. The following general result holds.

\begin{thm}\cite[Theorem 1.8]{Toma2}
\label{Toma2}
  Let $f(\vec{x})$ be a decomposable over $F_S$ form such that $f(\OO_S^m)$ is a discrete subset of $F_S$. Then $f(\vec{x}) = ag(\vec{x})$, where $g(\vec{x}) \in F[\vec{x}]$ and $a \in F_S^*$.
\end{thm}

It is easy to see that the analogous statement is not valid for non-decomposable forms. (The form $f(x_1,x_2,x_3) = (x_1^2 + \sqrt{2}x_2^2)x_3$ provides a simple counter-example.)
The proof of Theorem \ref{Toma2}  is based on the description of the closed $T$-orbits on $G/\Gamma$, where $G$ is the group of $F_S$-rational points $\mathbf{G}(F_S)$ of a reductive algebraic group $\G$ defined over $F$, $\Gamma$ is an arithmetic subgroup of $G$, and $T = \T(F_S)$, where $\T$ is a maximal $F$-split subtorus of  $\G$. (See \cite[Theorem 1.4]{Toma2} and \cite[Theorem 1.1]{TW}.)

In this paper the $S$-forms are arbitrary, \textit{not necessarily decomposable over }$F_S$.
Further on, $G = \underset{v \in S}{\prod}\mathrm{SL}_m(F_v)$, $G(F) = \mathrm{SL}_m(F)$, and $\Gamma = \mathrm{SL}_m(\mathcal{O}_S)$. Recall that $\Gamma$ is a lattice in $G$, the homogeneous space $G/\Gamma$ is endowed with the quotient topology, and
 $\pi: G \rightarrow G/\Gamma, g \mapsto g\Gamma$.
 The group $\mathrm{GL}_m(F)$ is acting naturally on $F_S[\vec{x}]$ and two forms $f$ and $f' \in F_S[\vec{x}]$ are called $F$-\textit{equivalent} if $f' = gf$, $g \in \mathrm{GL}_m(F)$. It follows from the definition of $f$ that $m \geq n$. If $(gf)(x_1, \dots, x_m) = (gf)(x_1, \dots, x_n, 0 \cdots, 0)$, $g \in \mathrm{GL}_m(F)$, then $f$ is $F$-\textit{equivalent to a form in $n$ variables}.

The next theorem reduces the study of forms in $m$ variables of degree $n$ to the case when $m = n$.
\begin{thm}
\label{reduction thm}
  If $f$ is an $S$-form of degree $n$ in $m$ variables and $f(\OO_S^m)$ is discrete in
  $F_S$ then $f$ is $F$-equivalent to a form in $n$ variables.
\end{thm}

\textit{Up to the end of the Introduction we assume that $m = n$.}
In order to formulate our main result, we recall the notion of $S$-\textit{norm form} and introduce the notion of $S$-\textit{quasi-norm form}. (If $F = \Q$ and $S$ is a singleton the notion of quasi-norm form is introduced in \cite{Toma4}.)
Let $K = F \alpha_1 + \dots + F \alpha_n$ be a field extension of degree $n$ of $F$ and $\{\theta_1, \dots, \theta_n\}$
be the set of all embeddings of $K$ into $\C$ over $F$, that is, $\theta_i(a) = a$ for $a \in F$.
Denote
$$l_i(\vec{x}) = \theta_i(\alpha_1)x_1 + \dots + \theta_i(\alpha_n)x_n.$$
Let
$\mathfrak{N}_{K/F}(\vec{x}) = \overset{n}{\underset{i=1}{\prod}}l_{i}(\vec{x})$ be the usual algebraic norm considered as a polynomial in $F_{S}[\vec{x}]$.
By an \textit{$S$-norm form corresponding to $K/F$}, we mean an $S$-form $F$-equivalent to
$$
f(\vec{x}) = c\cdot \mathfrak{N}_{K/F}(\vec{x}), \  c \in F_S^{*}.
$$

Now, let $n = 2s$ and $K$ be \textit{a totally real number field of degree} $s$ \textit{containing} $F$.
By an $S$-\textit{quasi-norm form corresponding to $K/F$} we mean a form $F$-equivalent to $f = (f_v)_{v \in S}$, where every $f_{v}$ is defined as follows.
 Consider the norm form $\mathfrak{N}_{K/F}(\vec{y})$, where $\vec{y} = (y_1, \dots, y_s)$. It is well-known
that
$$
K \underset{F}{\otimes}F_v = K_{v,1} \oplus \cdots \oplus K_{v,r_v},
$$
where $K_{v,i}$ are field extensions of $F_v$. Moreover,
 \begin{equation}
\label{decomp}
\mathfrak{N}_{K/F}(\vec{y}) = \underset{1 \leq i \leq r_v}{\prod}\mathfrak{N}_{K_{v,i}/F_v}(\vec{y})
\end{equation}
is the decomposition of $\mathfrak{N}_{K/F}(\vec{y})$ in a product of $F_v$-irreducible polynomials. For every $1 \leq i \leq r_v$ fix a quadratic form $q_{v,i}$ in two variables with coefficients from $F_v$ which is $F_v$-anisotropic, that is, $q_{v,i}(\vec{a}) \neq 0$ for $\vec{a} \in F_v^2 \setminus \{\vec{0}\}$. Denote
$\vec{x} = (\vec{x}_1, \vec{x}_2)$, where $\vec{x}_1 = (x_1, \dots, x_s)$ and $\vec{x}_2 = (x_{s+1}, \dots, x_n)$. Then $f_{v}$ is defined by
\begin{equation}\label{quasi-form}
f_{v}(\vec{x}) = \underset{1 \leq i \leq r_v}{\prod}q_{v,i}(\mathfrak{N}_{K_{v,i}/F_v}(\vec{x}_1), \mathfrak{N}_{K_{v,i}/F_v}(\vec{x}_2)).
\end{equation}

Our main theorem is the following.

\begin{thm}
\label{thm++} Let $f$ be an $S$-form and $H$ be the Zariski identity component of the stabilizer of $f$ in $G$. The following conditions are equivalent:
\begin{enumerate}
  \item $f$ is either an $S$-norm form or an $S$-quasi-norm form;
  \item $f(\OO_S^n)$ is discrete and $f$ does not represent $0$ over $F$ non-trivially;
 \item $H\pi(e)$ is compact.
\end{enumerate}
\end{thm}

For $F = \Q$ and $S = \infty$ the theorem is proved in \cite[Theorem 1.2]{Toma4}.

It follows from our definitions of $S$-norm and $S$-quasi-norm forms that an $S$-norm form corresponding to $K/F$ is also an $S$-quasi-norm form if and only if $K$ is a $CM$-field and $F$ is contained in the maximal totally real subfield of $K$, cf.\cite[Remark 4.4]{Toma4}.

Theorem \ref{thm++} implies immediately.

\begin{cor}
\label{cor1} Let $f(\vec{x})$ be a form with the following properties:
\begin{enumerate}
  \item[(i)]
$f(\OO_S^n)$ is discrete in $F_S$;
\item[(ii)]
$f$ does not represent $0$ over $F$ non-trivially;
\item[(iii)] $n$ is odd or $F$ is not a totally real number field.
\end{enumerate}
Then $f$ is an $S$-norm form.
\end{cor}

The following generalization of the conjecture of Cassels and Swinnerton-Dyer is natural.

\textbf{Conjecture}. \textit{Let $F$ be a totally real number field and $f(\vec{x})$ be a decomposable over $F_S$ form such that
\begin{enumerate}
  \item[(i)] $n \geq 3$;
\item[(ii)]
$f$ does not represent $0$ over $F$ non-trivially;
\item[(iii)] $f(\OO_S^n)$ is separated from $0$, that is, there exists a neighbourhood $W$ of $0$ in $f(F_S^n)$ such that $f(\OO_S^n) \cap W = \{0\}$.
\end{enumerate}
Then $f$ is an $S$-norm form corresponding to a field extension $K/F$ of degree $n$.}

In view of the $S$-adic version of Mahler's compactness criterium (Theorem \ref{S-adic Mahler}) the above conjecture admits the following reformulation in homogeneous dynamical terms: if $F$ is a totally real number field and $H$ is a maximal diagonalizable over $F_S$ subgroup of $\mathrm{SL}_n(F_S)$, $n \geq 3$, then any  relatively compact orbit $H\pi(g)$ on $G/\Gamma$ ($= \mathrm{SL}_n(F_S)/\mathrm{SL}_n(\OO_S)$) is compact.

The following result shows that a natural generalization of the conjecture for forms non-decomposable over $F_S$ does not hold.

\begin{thm} \label{cor2}\cite[Theorem 1.4]{Toma4} Let $n \geq 3$, $F = \Q$ and $S$ be a singleton. Then if $H$ is a maximal diagonalizable over $\C$ but not over $\R$ subgroup of $\mathrm{SL}_n(\R)$, there exists an orbit ${H\pi(g)}$ which is relatively compact but not compact. Equivalently, if $n = 2s+t$, $s > 0$, then there exists a real form $f$ of degree $n$ with $s$ the number of pairs of complex-conjugate linear forms dividing $f$ such that $f$ is neither a norm nor a quasi-norm form although $f$ does not represent $0$ over $\Q$ non-trivially and $f(\Z^n)$ is separated from $0$.
\end{thm}

\section{Preliminaries} \label{Preliminaries}
\subsection{Normalized valuations and norms.}\label{2.1}
For $v \in S$, let $|\cdot|_v$ be the \textit{normalized
valuation of} $F_v$. Recall that if $F_v = \C$ then $|a|_v$ is the square of the absolute value of the complex number $a$. %of  value of if $v$ is archimedean and the
The $F_S$-module $F_S^n$ is identified with the direct product of vector spaces $\underset{v \in S}{\prod}F_v^n$.
Every vector space $F_v^n$ is endowed with a \textit{normalized norm} $\| \cdot \|_v$ defined as follows\footnote{In this and the next sub-section we use notions and terminology introduced in \cite[Sections 2 and 3]{Toma2} and \cite[Section 5]{KT}.}. If $F_v = \R$ or $v$ is non-archimedean then every norm of the vector space $F_v^n$ is normalized. (Recall that all norms of a finite dimensional vector space over $F_v$ are equivalent.) In the case $F_v = \C$ the normalized norm $\| \cdot \|_v$ is defined as the square of an usual (for example, the Euclidean) norm of $\C^n$.
For $\vec{w} = (\vec{w}_v)_{v \in S} \in F_S^n$, we define the \textit{normalized norm of} $\vec{w}$ by
$$\|\vec{w}\|_S = \underset{v \in S}{\max}\| \vec{w}_v \|_v.$$
As  in \cite[Sections 2 and 3]{Toma2} and \cite[Section 5]{KT}, the \textit{content of} $\vec{w}$ is defined by
$$
\mathrm{cont}(\vec{w}) = \underset{v \in S}{\prod} \|\vec{w}_v \|_v.
$$
Similarly, we define \textit{content} for $a = (a_v)_{v \in S} \in F_S$ by
$$
\mathrm{cont}(a) = \underset{v \in S}{\prod} | a_v |_v.
$$
Since $\underset{v \in S}{\prod}|\xi|_v = 1$ for every $\xi \in \OO_S$ \cite[Ch.2, Theorem 12.1]{Cassels}, we have
$$
\mathrm{cont}(\vec{w}) = \mathrm{cont}(\xi\vec{w}) \ \textrm{for} \ \textrm{all} \ \xi \in \OO_S^*.
$$
The next lemma is proved in \cite[Lemma 3.2]{Toma2} for $s = 1$
and its proof for $s > 1$ is virtually the same.

\begin{lem}
\label{lemma from DJM} Let $s \in \N$. There exists a real $\kappa > 1$ such that for every $\vec{w} \in F_S^n$ with $\mathrm{cont}(\vec{w}) \neq 0$ there exists
$\xi \in \OO_S^*$ such that
$$
\frac{\mathrm{cont}(\vec{w})^{1/n}}{\kappa} \leq \|\xi^s\vec{w}\|_S \leq \kappa \cdot\mathrm{cont}(\vec{w})^{1/n}.
$$
\end{lem}

\subsection{$S$-adic compactness criterium.}\label{2.2}

Denote by $\mathcal{B}(a)$ the open ball in $F_S^n$ of radius $a$ centered at $0$ with respect to  $\|\cdot\|_S$. By a pseudo-ball in $F_S^n$ of radius $a$ centered at $0$, we mean the set $\mathcal{P}(a) = \{\vec{w}\in F_S^n: \mathrm{cont}(\vec{w}) < a\}$.

Recall the following $S$-adic version of Mahler's compactness criterium\footnote{To the best of our knowledge, the $S$-adic version of Mahler's criterium was first proved in \cite[Theorem 5.12]{KT}. A more general result, including the global fields in positive characteristic, is proved in \cite[Theorem 1.1]{KST}}.

\begin{thm}
\label{S-adic Mahler}
  Given a subset $M \in G$, the following conditions are equivalent:
  \begin{enumerate}
    \item $\pi(M)$ is relatively compact in $G/\Gamma$;
    \item there exists $a > 0$ such that $g\OO_S^n \cap \mathcal{P}(a) = \{0\}$ for all $g \in M$;
    \item there exists $a > 0$ such that $g\OO_S^n \cap \mathcal{B}(a) = \{0\}$ for all $g \in M$.
  \end{enumerate}
\end{thm}

\subsection{Algebraic groups and tori.}\label{2.3}
The group $\mathbf{G} = \mathbf{SL}_m$  is considered as an algebraic group defined over $F$. Hence $G = \mathbf{G}(F_S) = \underset{v \in S}{\prod}G_v$, where $G_v = \mathbf{G}(F_v)$. On every $G_v$ we have Zariski topology induced by the Zariski topology on $\mathbf{G}$ and Hausdorff topology induced by the Hausdorff topology of $F_v$. By Zariski (respectively, Hausdorff) topology on $G$ we mean the product of the Zariski (respectively, Hausdorff) topologies on $G_v$, $v \in S$. In order to distinguish the two topologies, the topological notions connected with the first one will be used with the prefix "Zariski".

If $L = \mathbf{L}(P)$ where $\mathbf{L}$ is an algebraic group defined over a local field $P$ of characteristic $0$ the Zariski closed subgroups of
$L$  will be also called \textit{algebraic subgroups of $L$}.

By a torus of $G_v$ we mean a Zariski connected abelian algebraic subgroup of $G_v$ consisting of semisimple elements. If $T_v, v \in S$, are tori of $G_v$ then $T = \underset{v \in S}{\prod}T_v$ is a torus of $G$. A torus $T$ of $G$ is an $F$-\textit{torus} if $T = \mathbf{T}(F_S)$, where $\mathbf{T}$ is an algebraic torus of $\mathbf{G}$ defined over $F$.
It is well-known that every $T_v$ is an almost direct product of its maximal $F_v$-anisotropic (i.e., without non-trivial $F_v$-characters) subtorus $T_{v,a}$ and of its maximal $F_v$-split (or, diagonalizable over $F_v$) subtorus $T_{v,s}$ \cite[\S8.15]{Borel}.

The following connection between the maximal tori and the $F_v$-forms is easy to prove. (A proof for $F_v = \R$ is available in \cite[Proposition 2.1]{Toma4} and the proof for any $F_v$ is virtually the same.)
\begin{prop}
\label{prop.inv}
Let $F_v[\vec{x}]_n$ be the space of all homogeneous
polynomials  of degree $n$ in $n$ variables $\vec{x} = (x_1, \dots, x_n)$. The following holds:
\begin{enumerate}
\item[(a)] If $H_v$ is a maximal torus of $G_v$ and $F_v[\vec{x}]_n^{H_v}$ is the subspace of all $H_v$-invariant polynomials then
\begin{equation}
\label{i-form}
F_v[\vec{x}]_n^{H_v} = F_vf_v,
\end{equation}
where $f_v$ is an $F_v$-form.
\item[(b)] Conversely, the Zariski identity component of the stabilizer of every $F_v$-form in $F_v[\vec{x}]_n$ is a maximal torus of $G_v$.
\end{enumerate}
\end{prop}

If $H_v$ is an algebraic subgroup of $G_v$ then $\mathfrak{R}_u(H_v)$ denotes the unipotent radial of $H_v$, that is, $\mathfrak{R}_u(H_v)$ is the maximal Zariski connected unipotent normal subgroup of $H_v$. If $H = \underset{v \in S}{\prod}H_v$ then $\mathfrak{R}_u(H) = \underset{v \in S}{\prod}\mathfrak{R}_u(H_v)$ is \textrm{the unipotent radical of} $H$.

\section{A reduction result}\label{reduction result}

\subsection{On the stabilizer of a polynomial.}\label{2.4}
The following proposition is proved in the real case (that is, for $G = \mathrm{SL}_m(\R)$ and $\Gamma = \mathrm{SL}_m(\Z)$) in
 \cite[Proposition 3.1]{Toma4}.

\begin{prop}
\label{connection1} Let $\mathrm{p} \in F_S[\vec{x}]$, and let $\mathrm{St}_G(\mathrm{p})$ be the stabilizer of $\mathrm{p}$ in $G$. Suppose that $\mathrm{p}(\OO_S^m)$ is
discrete in $F_S$. Then $\mathrm{St}_G(\mathrm{p})\pi(e)$ is closed in
$G/\Gamma$.
\end{prop}
{\bf Proof.} Let $h_i \in \mathrm{St}_G(\mathrm{p})$ be a sequence such that $\lim_{i\rightarrow \infty} h_i\pi(e) = \pi(g)$. We need to prove that $\pi(g) \in \mathrm{St}_G(\mathrm{p})\pi(e)$. There exists a sequence $\tau_i \in G$ converging to $e$ such that
\begin{equation}\label{poly}
h_i\pi(e) = \tau_i \pi(g),
\end{equation}
 equivalently,
$$
h_i \gamma_i = \tau_i g,
$$
where $\gamma_i \in \Gamma$. Let $\vec{z} \in \OO_S^m$. Then
$$
\mathrm{p}(\gamma_i \vec{z}) = \mathrm{p}(\tau_i g\vec{z})
$$
and, therefore, $\mathrm{p}(\gamma_i\vec{z}) \underset{i \rightarrow \infty}{\rightarrow}\mathrm{p}(g\vec{z})$. %Since $\OO_S^m$ is %$\Gamma$-invariant and
Since $\mathrm{p}(\OO_S^m)$ is
discrete in $F_S$, there exists a $c(\vec{z}) \in \N$ such that
$$
\mathrm{p}(\tau_i g\vec{z}) = \mathrm{p}(g\vec{z})
$$
whenever $i > c(\vec{z})$.

For every $\phi = (\phi_v)_{v \in S} \in F_S[\vec{x}]$, denote $\deg \phi = \max \{\deg \phi_v: v \in S\}$.  Let $r = \deg \mathrm{p}$ and $V$ be the $F_S$-module of all  $\phi \in F_S[\vec{x}]$ with $\deg \phi \leq r$. Since $\OO_S^m$ is countable, we can write $\OO_S^m = \{\vec{z}_i: i \in \N\}$.
Let $V_k = \{\phi \in V: \phi(\vec{z}_i) = 0 \ \textrm{for} \ \textrm{all} \ i \leq k\}$.
Then  $V_k = \underset{v \in S}{\prod}V_{k,v}$, where $V_{k,v} = \{\phi_v \in F_v[\vec{x}]: \deg(\phi_v) \leq r  \ \textrm{and} \ \phi_v(\vec{z}_i) = 0 \ \textrm{for} \ \textrm{all} \ i \leq k\}$. Since the natural imbedding of $\OO_S^m$ into $F_v^m$, $v \in S$, is Zariski dense in $F_v^m$, there exists
$l$ such that $V_{l,v} = \{0\}$ for all $v \in S$. Equivalently, $V_l = \{0\}$.
Put $i_0 = \max \{c(\vec{z}_i): 1 \leq i \leq l\} + 1$. It follows from the definitions of $c(\vec{z})$ and $l$ that
$$
\mathrm{p}(\tau_{i_0} g\vec{x}) = \mathrm{p}(g\vec{x}),
$$
that is, $\tau_{i_0} \in \mathrm{St}_G(\mathrm{p})$. In view of (\ref{poly}),
$$
\pi(g) = \tau_{i_0}^{-1}h_{i_0}\pi(e) \in \mathrm{St}_G(\mathrm{p})\pi(e).
$$
\qed
\subsection{Proof of Theorem \ref{reduction thm}.}\label{2.5}
It follows from Proposition \ref{connection1} that $H\pi(e)$ is closed. For every $v \in S$, let $H_{v,u}$ be the subgroup of $H_v$ generated by the $1$-parameter unipotent
subgroups of $H_v$. It is easy to see that $H_{v,u}$ is a Zariski closed normal subgroup of $H_v$ and $H_{v}/H_{v,u}$ consists of semisimple elements. Put
$H_u = \underset{v \in S}{\prod}H_{v,u}$.
 According to
\cite[Theorem 1]{Toma1} there exists a Zariski connected $F$-subgroup $\mathbf{P}$ of $\mathbf{G}$ and a subgroup of finite index $P'$ in $\mathbf{P}(F_S)$ containing $H_u$ and contained in $H$ such that
$$
\overline{H_u\pi(e)} = P'\pi(e),
$$
$P' \cap \Gamma$ is a lattice in $P'$ and, moreover, $\mathbf{P}$ has the property: if $\mathbf{L}$ is a proper normal $F$-subgroup  of $\mathbf{P}$ then there exists $v \in S$ such that $(\mathbf{P}/\mathbf{L})(F_v)$ contains a unipotent element different from the identity. Hence $H_u = P'$ and $H_u \cap \Gamma$ is an $S$-arithmetic lattice in $H_u$. Since $\mathcal{R}_u(H) = \mathcal{R}_u(H_u)$ we get that $\mathcal{R}_u(H) \cap \Gamma$ is Zariski dense in $\mathcal{R}_u(H)$.
 Let
$$
W = \{\vec{w} \in F_S^m: h\vec{w} = \vec{w} \ \mathrm{for} \ \mathrm{all} \ h \in \mathcal{R}_u(H)\}.
$$
Then the $F$-vector space $W \cap F^m$ is Zariski dense in $W$ and, clearly, $\dim_F (W \cap F^m) = m-n$.

Note that ${H}_u = {A} \ltimes \mathcal{R}_u({H}_u)$, where $A$ is the group of $F_S$-rational points of an $F$-algebraic group $\mathbf{A}$ isomorphic over $F$ to $\mathbf{SL}_{m-n}$, and $W$ is $A$-invariant.
Let
$$
W_0 = \{\vec{w} \in F_S^m: s\vec{w} = \vec{w} \ \mathrm{for} \ \mathrm{all} \ s \in A)\}.
$$
Since $\mathbf{A}(F)$ is Zariski dense in ${A}$, the $F$-vector space $W_0 \cap F^m$ is Zariski dense in $W_0$, $W_0 \cap F^m \cong F^n$, and $$F^m = (W_0 \cap F^m) \oplus (W \cap F^m).$$
It is easy to see that $H$ acts trivially on $W$ and $W_0$ is $H$-invariant. Hence if $\{\vec{v}_1, \dots, \vec{v}_m\}$ is a basis of $F^m$ such that $\{\vec{v}_1, \dots, \vec{v}_n\}$ is a basis of $W_0 \cap F^m$ and $\{\vec{v}_{n+1}, \dots, \vec{v}_m\}$ is a basis of $W \cap F^m$ then with respect to this basis the form $f$ is $F$-equivalent to a form $f'$ such that
$$
f'(y_1, \dots, y_m) = f'(y_1, \dots, y_n, 0, \dots, 0).
$$
\qed

\section{Proof of Theorem \ref{thm++}}\label{main result}

Further on $f$ is an $S$-form of degree $n$ on $F_S^n$ and $H$ is the Zariski identity component of the stabilizer of $f$ in $G$. According to Proposition \ref{prop.inv}  $H = \underset{v \in S}{\prod}H_v$, where $H_v$ is a maximal torus in $G_v$. %It is easy to see that $f$ is determined by $H$ up to a multiplication by an element from %$F_S^*$.

\subsection{A compactness criterium.}\label{2.6}

The following compactness criterium represents the $S$-adic version of a similar result \cite[Proposition 4.1]{Toma4} when $F = \Q$ and $S$ is a singleton.

\begin{prop}
\label{auxiliary2}  $H\pi(e)$ is compact if and only if $f(\mathcal{O}_S^n)$ is discrete in $F_S$ and $f(\vec{w}) \neq 0$ for every $\vec{w} \in F^n \setminus \{\vec{0}\}$.
\end{prop}

The proof is based on the following

\begin{lem}
\label{auxiliary2+}  Let $P$ be a local field of characteristic $0$ and $T$ be the maximal algebraic torus of $\mathrm{SL}_n(P)$ fixing a form $\phi$ of degree $n$. Let $T_s$ be the maximal $P$-split sub-torus of $T$, $\| \cdot \|$ be a norm on the vector space $P^n$, and $| \cdot |$ be a valuation of the local field $P$.
\begin{enumerate}
  \item[(i)]
 There exists a constant $k > 1$ such that for all $\vec{w} \in P^n$ with $\phi(\vec{w}) \neq 0$ there exists  $t \in T_s$ such that
  \begin{equation}\label{asymp1}
  \frac{1}{k} |\phi(\vec{w})|^{1/n} \leq \| t\vec{w} \| \leq {k}|\phi(\vec{w})|^{1/n}.
  \end{equation}
 \item[(ii)] If $\phi(\vec{w}) = 0$ then for every $\varepsilon > 0$ there exists $t_{\varepsilon} \in T_s$ such that
  $$
  \|t_{\varepsilon}\vec{w}\| < \varepsilon.
  $$
  \end{enumerate}
\end{lem}
{\bf Proof.} If $P = \C$ then the lemma follows easily from the fact that $\phi$ is diagonalizable over $P$. Further on we suppose that $P \neq \C$.

Let $R$ be the subspace of the linear space
$\mathrm{M}_n(P)$ generated by $T$. By theorems of Wedderburn \cite[Theorem 1.4.4]{Herstein} and Wedderburn-Artin \cite[Theorem 2.1.6]{Herstein} $R = R_1 \oplus \cdots \oplus R_r$, where $R_i$ is a field extension of $P$ of degree $n_i$. It is clear that $n = n_1 + \cdots + n_r$. We have $P^n = V_1 \oplus \cdots \oplus V_r$, where $V_i$
are subspaces of $P^n$, $\dim_P V_i = n_i$, and $R_i$ is a maximal subfield of $\mathrm{End}(V_i)$. Note that if $a \in R_i \subset \mathrm{End}(V_i)$ then $\det(a) = \mathfrak{N}_{R_i/P}(a)$.
For every $i$ we fix a non-zero vector $\vec{v}_i \in V_i$ and consider the map $\omega_i: R_i \rightarrow V_i, x \mapsto x(\vec{v}_i)$. It is easy to see that $\omega_i$ is an isomorphism of $P$-vector spaces and, moreover, considering $R_i$ as an $R_i$-module, $\omega_i$ is an isomorphism of $R_i$-modules.
Further on, $R_i$ is identified with $V_i$ via $\omega_i$.

Let $\vec{w} = \vec{w}_1 + \dots + \vec{w}_r \in P^n, \vec{w}_i \in V_i$, and $\phi(\vec{w}) \neq 0$. In view of Proposition \ref{prop.inv}, replacing $\phi$ by its constant multiple, we get
\begin{equation}\label{alg.norm}
  \phi(\vec{w}) = \mathfrak{N}_{R_1/P}(\vec{w}_1) \dots \mathfrak{N}_{R_r/P}(\vec{w}_r).
\end{equation}
The valuation $| \cdot |$ on $P$ is uniquely extended to a norm $\|\cdot\|_i$ on $V_i$ (identified with the field $R_i$), as follows
\begin{equation}\label{alg.norm+}
\|\vec{w}_i\|_i = |\mathfrak{N}_{R_i/P}(\vec{w}_i)|^{\frac{1}{n_i}}.
\end{equation}
(See \cite[ch.2, Theorem 10.1]{Cassels}.) Since the norms on a finite dimensional vector space over a local field $P$ are all equivalent, if $r = 1$ then (\ref{alg.norm}) and (\ref{alg.norm+}) imply $(\textrm{i})$ with $t$ equal to the identity map.

Suppose that $r > 1$. By the equivalence of the norms of the vector space $P^n$, we may (and will) suppose that the norm $\|\cdot\|$ on $P^n$ is defined by
\begin{equation}\label{norm on P^n}
\|\vec{w}\| = \underset{i}{\max}\|\vec{w}_i\|_i.
\end{equation}
Fix a $p \in P$ such that $|p| > 1$. Let $\mathcal{T} = \{(a_1, \dots ,a_r) \in \R_{+}^r: \underset{i = 1}{\overset{r}{\prod}}a_i^{n_i} = 1\}$ and let
$\Lambda = \{(\lambda_1, \dots, \lambda_r) \in \mathcal{T}: \textrm{every} \ \lambda_i \ \textrm{is} \ \textrm{a} \ \textrm{power} \ \textrm{of} \ |p|\}$. Then
$\Lambda$ is a co-compact lattice in the (multiplicative) Lie group $\mathcal{T}$. Hence there exists a constant $k > 1$ such that for every $(b_1, \dots ,b_r) \in \R_{+}^r$ there exists a $(\lambda_1, \dots, \lambda_r) \in \Lambda$ satisfying
\begin{equation}\label{asymp0}
  \frac{1}{k}\big(\underset{i = 1}{\overset{r}{\prod}}b_i^{n_i}\big)^{1/n} \leq \lambda_ib_i \leq {k}\big(\underset{i = 1}{\overset{r}{\prod}}b_i^{n_i}\big)^{1/n}.
\end{equation}
It follows from (\ref{alg.norm}) and (\ref{alg.norm+}) that
$$|\phi(\vec{w})| = \underset{i = 1}{\overset{r}{\prod}}\|\vec{w}_i\|_i^{n_i}.$$
Put $b_i = \|\vec{w}_i\|_i$. In view of (\ref{asymp0}), there exists $(l_1, \dots, l_r) \in \Z^r$ such that $n_1l_1 + \cdots + n_rl_r = 0$ and
\begin{equation}\label{asymp2}
  \frac{1}{k} |\phi(\vec{w})|^{1/n} \leq \| p^{l_i}\vec{w}_i \|_i \leq {k}|\phi(\vec{w})|^{1/n}.
\end{equation}
Let $t$ be the endomorphism of $P^n$ acting on every $V_i$ as a homothety with scale factor $p^{l_i}$. Then $t \in T_s$. Now, (\ref{asymp1}) follows immediately from  (\ref{asymp2}) and (\ref{norm on P^n}), proving $(\textrm{i})$.

Let $\phi(\vec{w}) = 0$. It follows from (\ref{alg.norm}) that $\vec{w}_{i_0} = 0$ for some $1 \leq i_0 \leq r$. Let $t' \in T_s$ be such that $t'$ acts on every $V_i, i \neq i_0,$ as a homothety with scale factor $d_i$ with $1 > d_i > 0$. It is clear that given $\varepsilon > 0$ we have that $|t_{\varepsilon}\vec{w}| < \varepsilon$ where $t_{\varepsilon} = t'^l$ with $l \in \N$ sufficiently large. This completes the proof of $(\textrm{ii})$. \qed

\medskip

In the formulation and the proof of the next corollary we will use the following conventional notation: if $\psi_1$ and $\psi_2$ are positive functions then the notation $\psi_1 \asymp \psi_2$ means that ${\psi_1}/{C} \leq \psi_2 \leq C\psi_1$ for some real constant $C > 1$. % which does not depend on the variables of $\psi_1$ and $\psi_2$.

\medskip

\begin{cor}
\label{auxiliary2++} Let $f = (f_v)_{v \in S}$ be an $S$-form with stabiliser $H = \underset{v \in S}{\prod}H_v$.
\begin{enumerate}
  \item[(i)] If $\vec{w} \in F_S^n$ and $\mathrm{cont}(f(\vec{w})) \neq 0$ then there exist $t \in H$ and $\xi \in \OO_S^*$ such that
  $$
  \mathrm{cont}(f(\vec{w}))^{1/n} \asymp \|t(\xi\vec{w})\|_S;
  $$
  \item[(ii)] If $\vec{w} \in F_S^n$ and $\mathrm{cont}(f(\vec{w})) = 0$ then for every $\varepsilon > 0$ there exist $t_\varepsilon \in H$ and $\xi_\varepsilon \in \OO_S^*$ such that
  $$
  \|t_\varepsilon(\xi_\varepsilon\vec{w})\|_S < \varepsilon.
  $$
\end{enumerate}
\end{cor}
{\bf Proof.} Let $\vec{w} = (\vec{w}_v)_{v \in S} \in F_S^n$ and $\mathrm{cont}(f(\vec{w})) \neq 0$. In view of Lemma \ref{lemma from DJM},
 there exists $\xi \in \OO_S^*$ such that
\begin{equation}\label{cont for f}
f_{v_1}(\xi\vec{w}_{v_1}) \asymp f_{v_2}(\xi\vec{w}_{v_2})
\end{equation}
for all
$v_1$ and $v_2 \in S$. By Lemma \ref{auxiliary2+}(i), there exists $t_v \in H_v$ such that
\begin{equation}\label{cont for f2}
|f_{v}(\xi\vec{w}_{v})|_v^{1/n} \asymp \|t_v(\xi\vec{w}_{v})\|_v
\end{equation}
for every $v \in S$. Put $t = (t_v)_{v \in S}$.
Now $(\textrm{i})$ follows from (\ref{cont for f}), (\ref{cont for f2}), and the equation $\mathrm{cont}(f(\xi\vec{w})) = \mathrm{cont}(f(\vec{w}))$.

Let $\mathrm{cont}(f(\vec{w})) = 0$. Then $f_{v_0}(\vec{w}_{v_0}) = 0$ for some $v_0 \in S$. Hence there exists a sequence $\xi_i \in \OO_S^*$ such that
$f_{v}(\xi_i\vec{w}_{v}) \rightarrow 0$ for every $v \in S \setminus \{v_0\}$. Let $\varepsilon_i \rightarrow 0^+$. Applying Lemma \ref{auxiliary2+}(i) when $f_{v}(\xi_i\vec{w}_{v}) \neq 0$, and Lemma \ref{auxiliary2+}(ii) when $f_{v}(\xi_i\vec{w}_{v}) = 0$, we get $t_{v,i} \in H_v$ such that $|f_{v}(\xi_i\vec{w}_{v})|_v^{1/n} \asymp \|t_{v,i}(\xi_i\vec{w}_{v})\|_v$ when $f_{v}(\xi_i\vec{w}_{v}) \neq 0$, and $\|t_{v,i}(\xi_i\vec{w}_{v})\|_v < \varepsilon_i$ when $f_{v}(\xi_i\vec{w}_{v}) = 0$. Put $t_i = (t_{v,i})_{v \in S}$. Then $t_i(\xi_i\vec{w}) \rightarrow 0$, which completes the proof of $(\textrm{ii})$. \qed

{\bf Proof of Proposition \ref{auxiliary2}.} Let $H\pi(e)$ be compact. Then
$$
H = \mathcal{H}(\Gamma \cap H),
$$
where $\mathcal{H}$ is a compact subset of $H$. Let $\vec{z} \in \OO_S^n \setminus \{\vec{0}\}$ and  $\mathrm{cont}(f(\vec{z})) = 0$. In view of Corollary \ref{auxiliary2++}$(\textrm{ii})$, there exist sequences $\xi_i \in \OO_S^*$ and $t_i \in H$ such that
$$
\|t_i(\xi_i\vec{z})\|_S \rightarrow 0.
$$
Let  $t_i = \tau_i \delta_i$, where $\tau_i \in \mathcal{H}$ and $\delta_i \in \Gamma \cap H$. Since $\mathcal{H}$ is compact, we may (and will) assume that
$$
\|\delta_i(\xi_i\vec{z})\|_S \rightarrow 0.
$$
The latter contradicts the discreteness of $\OO_S^n$ in $F_S^n$. Therefore $\mathrm{cont}(f(\vec{z}))$ $\neq 0$ for all $\vec{z} \neq \vec{0}$, in particular, $f$ does not represent $0$ over $F$ non-trivially.

Let $\{\vec{z}_i\}$ be a sequence in $\OO_S^n \setminus \{\vec{0}\}$ such that  $f(\vec{z}_i) \rightarrow a$, $a \in F_S$. Then $\mathrm{cont}(f(\vec{z}_i)) \rightarrow \mathrm{cont}(a)$. It follows from  Corollary \ref{auxiliary2++}$(\textrm{i})$ the existence of sequences $\xi_i \in \OO_S^*$ and $t_i \in H$ such that  $\|t_i(\xi_i\vec{z}_i)\|_S$ is a bounded sequence. As above, let $t_i = \tau_i \delta_i$, where $\tau_i \in \mathcal{H}$ and $\delta_i \in \Gamma \cap H$. Then the sequence
$\|\delta_i(\xi_i\vec{z}_i)\|_S$ is also bounded. Passing to a subsequence we may (and will) assume that $\delta_i(\xi_i\vec{z}_i) = \delta_{i_0}(\xi_{i_0}\vec{z}_{i_0})$ for all $i \geq i_0$. Therefore
$$
f(\vec{z}_i) = \xi_i^{-n}\xi_{i_0}^{n}f(\vec{z}_{i_0}),
$$
for all $i \geq i_0$. Since $f(\vec{z}_{i_0}) \neq 0$ and  $\OO_S^*$ is discrete in $F_S^*$ we get that $f(\vec{z}_i) = a$ for sufficiently large $i$. Therefore $f(\OO_S^*)$ is discrete in $F_S$.

Conversely, assume that $f(\OO_S^*)$ is discrete and $f$ does not represent $0$ over $F$ non-trivially. Then the orbit $H\pi(e)$ is closed by
Proposition \ref{connection1}. In view of Theorem \ref{S-adic Mahler}, it remains to prove that the set of lattices $\{h\OO_S^n : h \in H \}$ is separated from $0$. Suppose to the contrary that there exist $h_i \in H$ and $\vec{z}_i \in \OO_S^n \setminus \{\vec{0}\}$ such that $h_i (\vec{z}_i) \rightarrow \vec{0}$. Then
$f(\vec{z}_i) = 0$ for large $i$, which is a contradiction. \qed

\subsection{Completion of the proof of Theorem \ref{thm++}.}\label{2.7}
Theorem \ref{thm++} follows immediately from Proposition \ref{auxiliary2} and the next proposition.
\begin{prop}
\label{auxiliary!!!}  $H\pi(e)$ is compact if and only if $f$ is either an $S$-norm or an $S$-quasi-norm form.
\end{prop}
{\bf Proof.} %$\Rightarrow)$
Suppose that $H\pi(e)$ is compact. Let $\Delta$ be a subgroup of finite index without torsion of $H \cap \Gamma$. Denote by $K$ the subspace of the vector space $\mathrm{M}_n(F)$ generated by $\Delta$. It follows from theorems of Wedderburn \cite[Theorem 1.4.4]{Herstein} and Wedderburn-Artin  \cite[Theorem 2.1.6]{Herstein} that
$$K = K_1\oplus \cdots  \oplus K_r,$$
where $K_i$ are fields containing $F$. Let  $K^1 = K \cap H$, $\mathbf{T}$ be the identity component of the Zariski closure of $K^1$ in the $F$-algebraic group $\mathbf{G}$. Then $T = \mathbf{T}(F_S) = \underset{v \in S}{\prod}T_v \subset H$, where $T_v = \mathbf{T}(F_v)$.
Since $\Delta$ is a co-compact lattice in $H$, $H/T$ is compact. It follows from the compactness of $T/\Delta$ that $\mathbf{T}$ is an $F$-anisotropic algebraic torus, that is, $\mathbf{T}$ does not admit non-trivial $F$-characters. %It is easy to see that $\mathbf{T}$ is an $F$-anisotropic algebraic torus, that is, $\mathbf{T}$ does not admit
Note that, if $r \geq 2$ then the homomorphism $K^1 \rightarrow F^*, (a_1, \dots, a_r) \mapsto \mathfrak{N}_{K_1/F}(a_1)$ can be extended to a non-trivial $F$-character of $\mathbf{T}$. Therefore $r = 1$, that is,
$K$ is a subfield of $\mathrm{M}_n(F)$ containing its center $F$. Let $d = [K:F]$. It is well known that any maximal subfield of $\mathrm{M}_n(F)$ is an extension of $F$ of degree $n$. Hence, $n = dl$, where $l \in \N$. If $l = 1$, then $K$ is a maximal subfield of $\mathrm{M}_n(F)$, which implies that $T = H$ and $f$ is an $S$-norm form corresponding to the field extension $K/F$.

Let $l \geq 2$. In view of a theorem of Noether-Skolem \cite[Theorem 4.3.1]{Herstein}, $K$ is conjugated by an element from $\mathrm{GL}_n(F)$ to any subfield of $\mathrm{M}_n(F)$ isomorphic to $K$ and containing the center $F$ of $\mathrm{M}_n(F)$. Hence we may (and will) assume that the field $K$ is imbedded \textit{diagonally} in $\mathrm{M}_n(F)$ by a morphism $\mathfrak{m}: K \rightarrow \mathrm{M}_n(F)$ explicitly defined as follows. We fix an injective morphism $K \rightarrow \mathrm{M}_d(F), a \mapsto [a]$, and for every $a \in K$, we let %Now, every $a \in K$ corresponds to the matrix $\mathfrak{m}(a) \in \mathrm{M}_n(F)$ defined by
\begin{equation}\label{*}
\mathfrak{m}(a) =
\begin{pmatrix}
[a] & 0 & \cdots & 0 \\
0 & [a] & \cdots & 0 \\
\vdots  & \vdots  & \ddots & \vdots  \\
0 & 0 & \cdots & [a]
\end{pmatrix},
\end{equation}
where the number of $[a]$ in (\ref{*}) is equal to $l$.

Let $v \in S$. Recall that
$$
K\underset{F}{\otimes}{F_v} = K_{v,1} \oplus \cdots \oplus K_{v,r_v},
$$
where $K_{v,i}$ are the completions of the field $K$ with respect to its places over $v$.
If $[K_{v,i}:F_v] = d_{v,i}$ then $K_{v,i}$ is considered as a maximal subfield of $\mathrm{M}_{d_{v,i}}(F_v)$, $\underset{i = 1}{\overset{r_v}{\sum}}d_{v,i} = d$, and for every $a = (a_1, \dots, a_{r_{v}}) \in K\underset{F}{\otimes}{F_v}$, $a_i \in K_{v,i}$, we put
\begin{equation}\label{*+}
[a] =
\begin{pmatrix}
a_1 & 0 & \cdots & 0 \\
0 & a_2 & \cdots & 0 \\
\vdots  & \vdots  & \ddots & \vdots  \\
0 & 0 & \cdots & a_{r_{v}}
\end{pmatrix} \in \mathrm{M}_{d}(F_v).
\end{equation}
In edition, the notation (\ref{*}) with $a \in K$ will be also used with $a \in  K\underset{F}{\otimes}{F_v}$.

It follows from (\ref{decomp}) that
\begin{equation}\label{**}
T_v = \{\mathfrak{m}(a):\det [a] = \underset{i=1}{\overset{r_v}{\prod}}\mathfrak{N}_{K_{v,i}/F_v}(a_i) =  1\}.
\end{equation}
In view of (\ref{*}) and (\ref{*+}) the centralizer of
$T_v$ in $\mathrm{M}_n(F_v)$ is equal to a sub-algebra of $\mathrm{M}_n(F_v)$ isomorphic to $\underset{i=1}{\overset{r_v}{\oplus}}\mathrm{M}_{l}(K_{v,i})$. Since $H_v$ commutes elementwise with $K\underset{F}{\otimes}{F_v}$, we have
\begin{equation}\label{***}
H_v \subset G_v \cap \underset{i=1}{\overset{r_v}{\prod}}\mathrm{GL}_{l}(K_{v,i}).
\end{equation}

Let $H_{v,i} = H_v \cap \mathrm{GL}_{l}(K_{v,i})$ and $T_{v,i} = T_v \cap \mathrm{GL}_{l}(K_{v,i})$. Then $T_{v,i}$ is contained in the center of $ \mathrm{GL}_{l}(K_{v,i})$ and $H_{v,i}$ is a maximal sub-torus of $G_v \cap \mathrm{GL}_{l}(K_{v,i})$.  Moreover, $H_{v,i}/T_{v,i}$ is compact because $H_{v}/T_{v}$ is compact. Since $\mathrm{SL}_{l}(K_{v,i}) \subset G_v$ we get that $H_{v,i} \cap \mathrm{SL}_{l}(K_{v,i})$ is a maximal (algebraic) torus of $\mathrm{SL}_{l}(K_{v,i})$. But $T_{v,i} \cap \mathrm{SL}_{l}(K_{v,i})$ is finite. Therefore $H_{v,i} \cap \mathrm{SL}_{l}(K_{v,i})$ is a maximal compact torus of $\mathrm{SL}_{l}(K_{v,i})$.
Since $l \geq 2$ this is not possible when $K_{v,i} = \C$. Therefore $K$ is a totally real number field. Let $K_{v,i} = \R$. It is easy to see that $\mathrm{SL}_{l}(\R)$, $l \geq 2$, contains a compact maximal algebraic torus only if $l = 2$, in which case, the maximal \textrm{topological} torus $\mathrm{S0}_{2}(\R)$ of $\mathrm{SL}_{l}(\R)$ is also maximal \textrm{algebraic}. We have proved that $l = 2$.

It remains to show that the $S$-form $f = (f_v)_{v \in S}$ corresponding to $H$ is $F$-equivalent to a form for which $f_v$ is defined by (\ref{quasi-form}) for all $v \in S$. Since $l = 2$, in view of (\ref{*}), $F^n = V_1 \oplus V_2$, where $V_1$ and $V_2$ are subspaces of $F^n$ of dimension $\frac{n}{2}$ and they are isomorphic as $K$-modules. (Recall that by (\ref{*}) $K$ is a maximal subfield of $\mathrm{M}_{n/2}(F)$.) If $\vec{x}_1 = (x_1, \dots, x_{n/2})$ are the coordinates with respect to a basis of $V_1$ we let $\vec{x}_2 = (x_{n/2+1}, \dots, x_{n})$ be the coordinates with respect to the corresponding basis of $V_2$. Denote $\vec{x} = (\vec{x}_1, \vec{x}_2)$.

For every $1 \leq i \leq r_v$, we let
$$
\mathcal{P}_{v,i} = \{p \in F_v[\vec{x}] : p \ \mathrm{is} \ \mathrm{homogeneous} \ \mathrm{of} \ \mathrm{degree} \ d_{v,i}\} \cup \{0\}.
$$
The group $G_v$ is acting naturally on $F_v[\vec{x}]$ and, obviously, $\mathcal{P}_{v,i}$ is $H_v$-invariant. Recall that, in view of (\ref{decomp}), $\mathfrak{N}_{K_{v,i}/F_v}(\vec{y})$ is one of the irreducible factors of the decomposition of $\mathfrak{N}_{K/F}(\vec{y})$ over $F_v$. Let
$$
W_{v,i} = F_v \mathfrak{N}_{K_{v,i}/F_v}(\vec{x}_1) + F_v \mathfrak{N}_{K_{v,i}/F_v}(\vec{x}_2).
$$
Then $W_{v,i}$ is a $2$-dimensional subspace of $\mathcal{P}_{v,i}$ and
it follows from (\ref{*+}) that for all $\mathfrak{m}(a) \in T_v$ and $p \in W_{v,i}$, we have
\begin{equation}\label{weight}
\mathfrak{m}(a)p = \mathfrak{N}_{K_{v,i}/F_v}(a_i)p.
\end{equation}
Using (\ref{*}) (with $l = 2$), we also get
$$
W_{v,i} = \{p \in \mathcal{P}_{v,i} : \mathfrak{m}(a)p = \mathfrak{N}_{K_{v,i}/F_v}(a_i)p \ \mathrm{for} \ \mathrm{all} \  \mathfrak{m}(a) \in T_v\}.
$$
Hence $W_{v,i}$ is a weight space for the $F_v$-rational representation of $T_v$ on $\mathcal{P}_{v,i}$. This implies that $W_{v,i}$ is $H_v$-invariant. We get an
  $F_v$-rational homomorphism $\rho: H_v \rightarrow \mathrm{GL}(W_{v,i})$. Let $H_{v,a}$ be the maximal anisotropic sub-torus of the $F_v$-torus $H_v$. Then $H_v = H_{v,a} T_v$. Since $H_{v,a}$ is anisotropic and $\rho$ is rational, $\rho(H_v)$ is contained in an $F_v$-anisotropic torus of $\mathrm{GL}(W_{v,i})$. Therefore there exists a quadratic form $q_{v,i}$ in two variables with coefficients from $F_v$ and not representing $0$ over $F_v$ non-trivially such that $q_{v,i}(\mathfrak{N}_{K_{v,i}/F_v}(\vec{x}_1), \mathfrak{N}_{K_{v,i}/F_v}(\vec{x}_2))$ is $H_{v,a}$-invariant. This implies that the form $f_v$ defined by (\ref{quasi-form}) is $H_{v,a}$-invariant. In view of (\ref{**}) and (\ref{weight}) $f_v$ is also $T_v$-invariant, consequently, $f_v$ is $H_v$-invariant. Therefore the form $f = (f_v)_{v \in S}$ is fixed by $H$ and it is an $S$-quasi-norm form corresponding to the field extension $K/F$.

Suppose that $f$ is a norm or a quasi-norm form corresponding to a field extension $K/F$. In the case of a norm form, $f = a\cdot \mathfrak{N}_{K/F}$. Then $f$ does not represent $0$ nontrivially over $F$ and $f(\OO_F^n)$ is discrete in $F_S^n$ because $\mathfrak{N}_{K/F}(\vec{x}) \in F[\vec{x}]$. Hence $H\pi(e)$ is compact by Proposition \ref{auxiliary2}. Suppose that
$f = (f_v)_{v \in S}$ is a quasi-norm form corresponding to a field extension $K/F$ of degree $\frac{n}{2}$, where every $f_v$ is defined by (\ref{quasi-form}). Then $f$ is $K^1$-invariant where $K^1$ is the group of $F$-rational points of an anisotropic over $F$ sub-torus $\mathbf{T}$ of $\G$.
Hence if $T = \T(F_S)$ then $f$ is $T$-invariant and $T\pi(e)$ is compact. Since the quadratic forms $q_{v,i}$ does not represent $0$ over $F_v$ we get that
 $H/T$ is compact which implies the compactness of $H\pi(e)$. \qed

\end{document}